\documentclass[multhd]{svmult}

\usepackage{graphicx}         
\usepackage[bottom]{footmisc} 
\usepackage{amssymb}

\begin{document}

\title*{Sequential and Asynchronous Processes \\
Driven by Stochastic or Quantum Grammars\\
and their Application to Genomics: a Survey}
\titlerunning{Stochastic or quantum grammar-driven processes }

\author{Dimitri Petritis}
\institute{
Institut de recherche math\'{e}matique de Rennes\\
Universit\'{e} de Rennes 1 and CNRS, France\\
\texttt{dimitri.petritis@univ-rennes1.fr}
}

\maketitle

\begin{abstract}
We present the formalism of sequential and asynchronous processes
defined in terms of random or quantum grammars and argue that
these processes have relevance in genomics. To make
the article accessible to the non-mathematicians, we keep
the mathematical exposition as elementary as possible, focusing
on some general ideas behind the formalism and stating the implications
of the known mathematical results. We close with a set of open
challenging problems.

\vskip3mm
\noindent
\textsl{Presented
at the European Conference of Mathematical and Theoretical
Biology, Dresden 18--22 July 2005.}


\end{abstract}


\section{The Classical Combinatorial Description: Configurations 
and Observables}
\label{sec:comb}
Biological molecules that convey information (nucleic acids and proteins)
can be abstractly modelled as finite sequences of letters chosen from an
alphabet $\mathbb{A}$. This alphabet can be of 4 letters, when we deal
with nucleic acids, of 20 letters, when we deal with proteins, or some
extended version of those alphabets --- possibly countable --- when some
additional information is incorporated into the description (eg. assigning 
exon-intron character to each nucleotide, including stereo-chemical
and pairing information to each amino-acid, etc.) 
All living organisms use the same alphabet to encode the biological information
relevant to their survival both as individuals and as species. Individuals
of the same species use sequences of (almost) fixed length.

To capture the combinatorial variety of all possible sequences arising
in all living matter,
we
use a universal set, the
\emph{(sequence) configuration space}.
An alphabet $\mathbb{A}$ 
gives rise to a configuration space $\mathbb{A}^*=\cup_{n=0}^\infty
\mathbb{A}^n$. 
Every possible and imaginable sequence is bijectively
mapped to a single point of $\mathbb{A}^*$. The configuration space
has a rooted tree structure carrying thus several natural distances:
the \emph{tree distance} counts the number of generations
one must go back to find a common prefix, 
the \emph{Hamming distance} counts the number of sites where
the residues differ.

An \emph{observable} $X$, taking values in some space $\mathbb{X}$, is
a map $X: \mathbb{A}^*\rightarrow \mathbb{X}$. In most situations,
the set $\mathbb{X}$ is just $\mathbb{R}$ or some subset of it. The most
elementary observables are $\{0, 1\}$-valued observables also called 
\emph{questions}, i.e.\ indicator
functions of subsets of $\mathbb{A}^*$.

Having delimited the playground --- the configuration space --- on which
all sequences can be represented, chronological changes occurring 
on sequences can be represented as a \emph{time evolution} on $\mathbb{A}^*$.
Here the term evolution 
applies mathematically to every time
scale to denote a $\mathbb{A}^*$-valued continuous time process, including
\begin{enumerate}
\item the accretion processes consisting in duplication of
the nucleic acid molecule, nucleotide by nucleotide, occurring
during cell division viewed as a 
very rapid length increasing time evolution
over nearest neighbouring points of the configuration space;
\item
the local modifications of the genetic sequence occurring
at random epochs of the life span of an individual such as
point mutations (DNA methylation, forward or backward 
slippage, etc.), modifying only few adjacent letters, viewed as
rapid evolution over closely lying points of
the configuration space 
(almost) preserving  the length of sequences;
\item
the global shuffling of vast regions of the genetic sequence 
of a given species occurring during trans-locations, inversions, independent
assortments or chiasmata during sexual reproduction viewed as rapid
length (almost) preserving time evolutions over distant points of the
configuration space but lying at (almost) the same distance from the root;
\item
the modifications of species sequences occurring over geological time
scales viewed as length non preserving slow evolution over distant
points of the configuration space.
\end{enumerate}
Since the configuration space $\mathbb{A}^*$ is countable, time evolution
is a continuous time jump process over $\mathbb{A}^*$; when sampled at the
instants of occurrence of the jumps, this process is 
a random walk on the tree $\mathbb{A}^*$.
For definiteness, we
focus only on evolution of sequences of nucleotides and more precisely
of the two first local types in the above list.

\begin{example}
Let $\mathbb{A}=\{A,C,G,T\}$. Then $\mathbb{A}^*=\cup_{N=0}^\infty
\mathbb{A}^N$, where 
$\mathbb{A}^0$ contains the empty sequence (denoted $\kappa$ in the sequel),
$\mathbb{A}^1$ contains the 4  sequences of length 1: $A,C, G$, and $T$, etc.
The set $\mathbb{A}^N$ contains $4^N$ sequences
of length $N$. The sequences of given
length can be represented as vertices of a given generation of a rooted
tree; all vertices but the root  have degree $|\mathbb{A}|+1=5$.
(See Fig.~\ref{fig:tree}.)
\begin{figure}
\centerline{\hbox{
\includegraphics[width=0.5\textwidth]{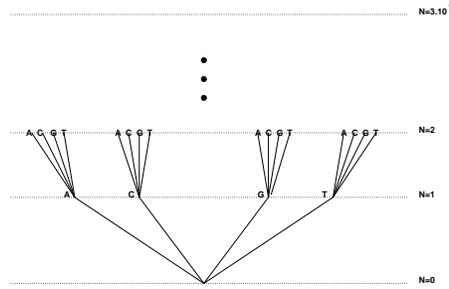}}}
\caption{\label{fig:tree} The bijection between the sequence configuration
space $\mathbb{A}^*$ and a rooted tree. To the root corresponds
the empty sequence $\kappa$ of zero length and no letter. 
To get the sequence mapped to another
vertex of the tree, one has to read the letters encountered on the path
joining the root to this particular vertex.}
\end{figure}
Sequences $AG$ and $GA$ have  tree and Hamming distance 2;
$AA$ and $GA$ have tree distance 2, Hamming distance 1.
\end{example}

\section{The Dual Probabilistic Classical  Description: States}
The configuration space as playground where sequences can be described
and evolve is very useful, concise, and powerful. Nevertheless, it is
much too precise to be directly exploitable: giving a point
in $\mathbb{A}^*$ for a
sequence of $N$ residues 
is equivalent to determining one out of $|\mathbb{A}|^N$ variables; for human
genome for instance the latter reads approximately $4^{3\times 10^9}$, a
tremendously huge number. Now when examining the DNA sequences appearing
in the cells of a multi-cell individual, we find that they are identical
but for some accidental modifications. 
What is needed is some qualitative description
of sequences allowing to make statistical comparisons between populations
(i.e.\ sets) of sequences. 
The configuration space is also
 too vast: the overwhelming majority of elements
of $\mathbb{A}^*$ never occur as  biologically viable and relevant sequences.
Finally,  the  evolution process 
on $\mathbb{A}^*$, even when it starts from a single sequence $\alpha \in 
\mathbb{A}^*$, will introduce some natural randomness on the set
of sequences obtained by the evolution of $\alpha$. 

The natural mathematical way to overcome all these difficulties is to deal
with probability measures, also called
\emph{states}, on  sets of sequences. A population of individuals instead
of being described by a precise set of elements of $\mathbb{A}^*$
is described by a state $\rho$, supported by this particular set, containing
all useful and relevant statistical information about
individual sequences.

Therefore, the pair consisting of a classical configuration space
and  a state  $(\mathbb{A}^*, \rho)$
is  nothing else than a
standard probability space; classical  observables
are nothing else than random variables on this space.

\section{On the Relevance of Quantum Mechanics in Biology: what, how, why?}

\emph{What?}
Mathematically, quantum mechanics is an extension of classical probability
theory where real random variables are replaced by self-adjoint operators
over an appropriate Hilbert space 
and classical states 
by quantum states.

Physically, quantum mechanics is a refinement of classical mechanics containing
a fundamental physical constant, the Planck constant 
$\hbar=1.05457\times 10^{-34}$Js. When action 
 values are comparable to  $\hbar$, the laws
of classical mechanics do not any longer describe reality satisfactorily.
Chemical properties like affinities of different atoms to form molecules, 
stability of matter,
conformational geometry of molecules, etc.\
rely on quantum mechanical rules. 

\emph{How?}
Biology relies on chemistry. Therefore it is pointless to ask
whether biology relies on quantum mechanics, it certainly does. The relevant
question is how it can. Due to the smallness of
$\hbar$, quantum phenomena become in general irrelevant when considering
large systems at high temperatures, and biomolecules are large warm systems!
Such systems although microscopically governed by quantum laws are globally
in the realm of classical physics. 
This phenomenon is known as \emph{decoherence}.
High temperature superconductivity \cite{BednorzMueller}
and
double-slit interference patterns for larger and larger molecules like
fullerenes $C_{60}$
\cite{ArndtNairzKellerZeilinger},
fluorofullerenes $C_{60}F_{48}$,
and even
biomolecules (meta-tetraphenylporphyrin $C_{44}H_{30}N_4$ involved
in hemoglobin transfers) observed lately, demonstrate coherent mesoscopic
quantum behaviour.

Several authors advocated that similar
emerging mesoscopic coherent quantum behaviour can arise 
in some circumstances in biology 
\cite{AlKhaliliMcFadden,Patel,Penrose}.
It is therefore worth studying biomolecules as quantum mesoscopic
systems. Notice however that we don't claim to solve
Schr\"odinger equation for every internal degree of freedom
of a given biomolecule, 
but rather study
the implication of emerging mesoscopic quantum behaviour with 
very few quantum degrees of freedom.

\emph{Why?}
All quantum phenomena at the size of a biomolecule in vivo, if any,
are necessarily fragile emergent ones, very easily returning to the classical
realm. It is estimated that for some biomolecules coherence can persist
up to  $10^{-5}$s; in \cite{AlKhaliliMcFadden} this estimate
is used to explain the rapidity of adaptive mutation occurring
in Escherichia coli in lactose environment.
It is in fact known that the time needed for the combinatorial
exploration of the configuration space is exponential in the length
of the searched sequence for classical search but only polynomial
for quantum search \cite{Kempe,HedaPetritis}. Moreover,
hitting time can be described as a partial measurement process
in quantum mechanics. Hence, quantum evolution and 
quantum measurement, 
provide  really new insight and explanation's for biological phenomena, 
strictly inside the known laws of Nature:
quantum evolution explains the 
observed rapidity of phenomena like adaptive mutation or emergence of life;
quantum measurement explains the process of selection of fittest or 
the adaptation  to the environment.

\section{Quantum Direct and  Dual Descriptions: Rays and Quantum States}
We give here
an elementary and very concise introduction to quantum description.
The interested reader can consult any standard book
on the subject (for instance \cite{Varadarajan}) or the
the freely available lecture notes \cite{Pet-QI} on the home page of the 
author. 

Quantum configuration spaces are complex separable Hilbert 
spaces.
For our purposes, the
\emph{quantum sequence configuration space} is the Hilbert space
$\mathbb{H}=\ell ^2(\mathbb{A}^*)$. An element of 
$\mathbb{H}$ is a complex function $\psi$ on  $\mathbb{A}^*$
such that
$\sum_{\alpha\in \mathbb{A}^*} |\psi(\alpha)|^2 <\infty$; the
scalar product is defined for all $\phi, \psi\in \mathbb{H}$ by
$\langle\ \phi\ |\ \psi\ \rangle=\sum_{\alpha\in \mathbb{A}^*} \overline{\phi}(\alpha)
\psi(\alpha)$.
An orthonormal basis of this space is provided by a family
$(e_\alpha)_{\alpha\in \mathbb{A}^*}$ of elements of $\mathbb{H}$, indexed
by the classical configurations. Therefore, basis elements are 
functions $ e_\alpha$ defined by $e_\alpha(\beta)=\delta_{\alpha, \beta}$,
for every $\alpha, \beta\in \mathbb{A}^*$, where $\delta$ denotes the Kronecker
symbol; this basis is isomorphic to
 $\mathbb{A}^*$. Every element  $\psi\in \mathbb{H}$ can be decomposed
as $\psi=\sum_{\alpha\in \mathbb{A}^*} \psi(\alpha) e_\alpha$.
Quantum configurations are \emph{rays}, i.e.\
vectors $\psi\in\mathbb{H}$ 
of unit norm.

Evolution is a unitary operator $U$ acting on $\mathbb{H}$, i.e.\ 
verifying $U^*U=UU^*=1$. Due to
the linear structure of the space, it is enough to study the action
of $U$ on the basis elements: it becomes then an infinite
unitary matrix. 

\emph{Quantum observables} are self-adjoint operators 
$X$ acting on $\mathbb{H}$, i.e.\ verifying $X^*=X$. 
On the basis vectors, they are represented
by infinite self-adjoint matrices.
The most elementary observables are projections (the quantum
analogue of indicators);
the spectral theorem establishes the decomposition of any self-adjoint
operator in terms of projections. 

\emph{Quantum states}, $\rho$, are self-adjoint, 
positive,
trace class,
normalised
operators acting on $\mathbb{H}$, called 
\emph{density matrices}.
Quantum observables are merely non-commutative random variables in the sense
that the expectation of $X$ in state $\rho$ is given by $\mathbb{E} X=
\mathop{\mathsf{tr}}(\rho X)$. Classical probability theory is a special
case of  quantum mechanics where all observables and states are represented
by diagonal matrices.

The measurement process is what renders quantum mechanics so counter-intuitive. 
Measurement corresponds to asking a question on the values an observable
can take. 
Suppose that we consider an observable $X$ taking a
discrete set of values $(\lambda_i)_i$. Therefore, we have
$X=\sum_i \lambda_i P_i$, where classically $P_i=1_{\{X=\lambda_i\}}$ 
 while quantum mechanically $P_i$ is the projector to the eigenspace
corresponding to the eigenvalue $\lambda_i$. Now perform the measurement
in the classical state (probability measure) 
or the quantum state (density matrix)
$\rho$ and consider the state after measurement in the two following 
situations: 
\begin{enumerate}
\item The result is filtred to get a conditioning on the fact 
that a particular value $\lambda_j$ has been observed 
after measurement: 
classically we get\\
 $\rho(\cdot)\mapsto  \rho(\cdot|X=\lambda_j)$, quantum
mechanically 
$\rho\mapsto \frac{P_j\rho P_j}{\mathop{\mathsf{tr}}(\rho P_j)}$.
\item The result is not filtred: classically we get, by virtue of Bayes' formula\\
 $ \rho(\cdot)\mapsto \sum_j\rho(\cdot|X=\lambda_j)\rho(X=\lambda_j)
=\rho(\cdot)$,\\ quantum mechanically
$\rho\mapsto \sum_j\frac{P_j\rho P_j}{\mathop{\mathsf{tr}}(\rho P_j)}
\mathop{\mathsf{tr}}(\rho P_j)=\sum_j P_j \rho P_j$.
\end{enumerate}
If $\rho$ and $X$ are not 
simultaneously diagonalisable, 
the state after measurement 
without filtering   $\sum_j P_j \rho P_j\ne\rho$.

\section{Grammars and Languages}
Grammars are powerful tools to drastically reduce the size of the
available configuration space by introducing 
combinatorial constraints by a set of elementary transformations
acting on the configurations (classical or quantum). 
In the sequel we use the symbol $\alpha$ to denote indistinguishably 
the classical
configuration $\alpha$ in the classical case and the
basis element $e_\alpha$ in the quantum case.

Grammars act on configurations in differents ways: sequentially in the case
of generational grammars,
in parallel for the so-called Lindemayer systems,
or in an asynchronous way. 

\subsection{Generational Grammars and their Classification}
Initially introduced for the description and analysis of natural
languages \cite{Chomsky}, generational grammars are extensively used nowadays
as models of computation.
A generational grammar $\Gamma$ is a small set of rules that
act sequentially on  elements of  $\mathbb{A}^*$  to produce a distinguished
subset of  $\mathbb{A}^*$, called the \emph{language} $L(\Gamma)$ generated
by $\Gamma$.
More specifically:
\begin{definition}
A \emph{(generational) grammar} $\Gamma$ is a quadruple $\Gamma=(\mathbb{A}_n, 
\mathbb{A}_t, \Pi, S)$, where
\begin{enumerate}
\item
$\mathbb{A}_n$ and 
$\mathbb{A}_t$ are two disjoint finite sets, the \emph{alphabets}
of non-terminal and terminal symbols respectively; we denote
by $\mathbb{A}= \mathbb{A}_n \cup\mathbb{A}_t$ and to avoid trivialities
we always assume that $\mathbb{A}_n\ne\emptyset$,
\item
$\Pi\subseteq (\mathbb{A}^+\setminus \mathbb{A}_t^*)\times \mathbb{A}^*$
is a finite set, the \emph{productions}, and
\item
$S\in \mathbb{A}_n$ is the initial symbol or \emph{axiom}.
\end{enumerate}
\end{definition}
Productions are rules for possible substitutions of subwords
of a sequence by other subwords. They define a binary relation on 
$\mathbb{A}^*$ as follows:
\begin{definition}
Let $\alpha, \beta\in \mathbb{A}^*$. We say that $\beta$ is 
\emph{directly derivable} from $\alpha$, and denote by $\rightsquigarrow$, the
binary relation 
$(\alpha \rightsquigarrow \beta)
\Leftrightarrow (\exists \alpha'\sqsubseteq \alpha; \exists
\beta'\sqsubseteq \beta: (\alpha',\beta')\in \Pi),$
where $\alpha'\sqsubseteq \alpha$ means that $\alpha'$ is a subword
of $\alpha$.
\end{definition}
Denote by $\stackrel{+}{\rightsquigarrow}$ the transitive 
closure\footnote{The \emph{transitive closure} of a binary
relation $R$ is the relation 
$R^+=\cup_{n=1}^\infty R^n.$} 
and by  $\stackrel{*}{\rightsquigarrow}$ the accessibility 
relation\footnote{The \emph{accessibility relation} associated with a
binary relation $R$ on $A$ is the binary relation  $R^*$ defined
for $a,b\in A$ by:
$(a R^* b)\Leftrightarrow ((a=b) \vee (a R^+ b)).$} 
of the direct derivability relation. 
\begin{definition} Let $\Gamma$ be a grammar.
The \emph{language} $L(\Gamma)$ generated by $\Gamma$ is the set
$L(\Gamma)=\{\alpha \in \mathbb{A}^* : 
S \stackrel{*}{\rightsquigarrow}\alpha\}.$
\end{definition}
We denote by $\textsf{Dom}(\Pi)=\{\alpha\in \mathbb{A}^+\setminus
\mathbb{A}_t^* \mid \exists \beta \in \mathbb{A}^* : (\alpha,\beta)\in\Pi\}$ 
and
$\textsf{Ran}_\Pi(\alpha)=
\{\beta\in \mathbb{A}^* | (\alpha,\beta)\in\Pi\}$ if $ \alpha\in 
\textsf{Dom}(\Pi)$ and 
$\emptyset$ otherwise.

If for some $\alpha\in\mathsf{Dom}(\Pi)$, and some $n$ we have
$\textsf{Ran}_\Pi(\alpha)=\{\beta^{(1)},\ldots,\beta^{(n)}\}$, then we use
the Backus-Naur shorthand notation $\alpha\rightarrow 
\beta^{(1)}|\cdots|\beta^{(n)}$ to mean that 
$ (\alpha,\beta^{(1)}),\ldots,(\alpha,\beta^{(n)})$ are all the 
possible productions with first element $\alpha$.
Grammars are classified according to their descendance degree
and their Chomsky (acontextuality) degree. 

The \emph{descendance degree}, $d$, is defined by
$d=\max\{|\textsf{Ran}_\Pi(\alpha)|: \alpha\in \textsf{Dom}(\Pi)\}$. Notice
that for all non-trivial grammars $d\geq 1$. If
$d=1$ the grammar is called \emph{deterministic} otherwise it is called
\emph{non-deterministic}.
It is worth noticing that deterministic descendance means that there exists
a function $\phi: \textsf{Dom}(\Pi)\rightarrow \mathbb{A}^*$
whose graph is the set $\Pi$, i.e.\ $(\alpha,\beta)\in \Pi\Leftrightarrow
\beta= \phi(\alpha)$. 
Fixing a given unpdating policy, this function induces
a dynamical system $\Phi: \mathbb{A}^*\rightarrow \mathbb{A}^*$
such that the sequence occuring as 
successive transformations of the grammar on the axiom $S$ appear as
the trajectory of the dynamical system $\Phi(S), \Phi\circ\Phi(S),
\ldots$. Such a trajectory is called a \emph{computational path} and it
can be finite if the system halts or infinite if it never halts.
For non-deterministic descendance, there does not exist such a function $\phi$,
or more precisely, this function is multi-valued. At each step, we must
use a branch of this function.
The branches are assigned a probability
vector or a unitary probability amplitude vector; we speak then
of a \emph{stochastic} or \emph{quantum} grammar respectively.
More precisely, if $\alpha\in\textsf{Dom}(\Pi)$, 
 stochastic descendance means that with 
$\alpha$ is associated a vector
$\vec{p}_\alpha=(p_{\alpha,\beta}, \beta\in \mathbb{A}^*)$ such that
$p_{\alpha,\beta}\geq 0$, $p_{\alpha,\beta}=0$ if $(\alpha,\beta)\not\in\Pi$,  
and  $\sum_\beta  p_{\alpha,\beta}=1$;
quantum descendance means that 
 with 
$\alpha$ is associated a vector
$\vec{u}_\alpha=(u_{\alpha,\beta}, \beta\in \mathbb{A}^*)$ such that
$u_{\alpha,\beta}\in\mathbb{C}$, $u_{\alpha,\beta}=0$ if $(\alpha,\beta)\not\in\Pi$,  
and  $\sum_\beta  |u_{\alpha,\beta}|^2=1$.

The second classifying parameter of grammars is their \emph{acontextuality}
or Chomsky degree. The different types of acontextuality are described in the
Table~\ref{tab:chomsky}. 
\begin{table}
\begin{tabular}{|c|l|l|l|l|}
\hline
Chomsky 
& Grammar & All productions $(\alpha,\beta)$ of the form & Recognition\\
\hline
0& recursively enumerable & $\alpha\rightarrow \beta $
with $\alpha\in \mathbb{A}^+\setminus\mathbb{A}^*_t, \beta\in\mathbb{A}^*$ & eTM\\
1& context-sensitive   
& $\alpha=\alpha_1\alpha'\alpha_2$ with $\alpha'\in\mathbb{A}_n^1$, 
$\alpha_1\alpha_2\ne\kappa$, & \\
&&$\beta= \alpha_1\beta'\alpha_2$, $\beta'\ne\kappa$
& eLBA \\
2& context-free     
& $\alpha\in \mathbb{A}_t^1$, $\beta\in \mathbb{A}^*$
& ePDA\\
3& regular      
& $\alpha\in \mathbb{A}_t^1$, $\beta\in \mathbb{A}_t^*$ or
$\beta \in   \mathbb{A}_n^*\times \mathbb{A}_t^*$
& eFA\\
\hline 
\end{tabular}
\caption{\label{tab:chomsky} The Chomsky hierarchy of grammars.
For every degree of acontextuality a universal
automaton can be used to recognise the language:
Turing machines (TM), linear nounded automata (LBA), push down automata (PDA),
or finite automata (FM).
For every descendance type,
the corresponding automaton acquires a prefix $e\in\{D, N, S,  Q\}$ meaning that
the evolution is deterministic, non-determintic (combinatorial), stochastic, or
quantum}
\end{table}
\begin{example}
A stochastic context-free grammar has been used in \cite{Sakakibaraetal}
to describe the secondary structure of RNA  molecule.
Its alphabets are $\mathbb{A}_t=\{A, C, G, U\}$ and 
$\mathbb{A}_n=\{S_0, \ldots, S_{13}\}$, the initial symbol $S_0$; its 
productions $\Pi$ are of the form
\begin{center}
{\scriptsize
\begin{tabular}{|l|l|l|l|l|}
$S_0\rightarrow S_1$  & $S_1\rightarrow CS_2G|AS_2U$ &
$S_2\rightarrow AS_3U$ & $S_3\rightarrow S_4S_9$ &$S_4\rightarrow US_5A$\\
$S_5\rightarrow CS_6G$ & $S_6\rightarrow AS_7$& 
$S_7\rightarrow US_7|GS_8$& $ S_8\rightarrow G|U  $& $ S_9\rightarrow AS_{10}U  $\\
 $S_{10}\rightarrow CS_{10}G|GS_{11}C$ &
 $S_{11}\rightarrow AS_{12}U  $&
$S_{12}\rightarrow US_{13}   $&
 $ S_{13}\rightarrow C  $&
\end{tabular}
}
\end{center}
A probability vector is associated with every production.
Every computational path leads to a different
realisation of the secondary structure. The elementary probability vectors
of the productions induce a natural probability measure on the
set of all possible secondary structures.
 A particular random
realisation
gives rise to the secondary structure depicted in Fig.~\ref{fig:rna}. 
\begin{figure}
\sidecaption
\centerline{\hbox{
\includegraphics[width=0.3\textwidth]{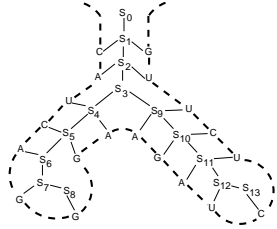}}}
\caption{\label{fig:rna} 
The secondary RNA structure 
as a particular random realisation derived by 
the generetional context-free stochastic 
grammar introduced in \cite{Sakakibaraetal}}
\end{figure}
\end{example}

Several other uses of accretion context-free stochastic grammars are reported
in bilogical literature.
Accretion dynamics defines
a random walk on the vertices of the computational paths leading
to the words of the language.
For stochastic context-sensitive grammars, the probability vectors
depend on the position of the random walk on the subtree of $\mathbb{A}^*$.
In this situation we speak about a  \emph{random environment}.
Several very specific models of random walks in random environment on trees
corresponding to particular classes of
context-sensitive grammars  have been studied in the literature
(the interested reader can look at  \cite{ComMenPop,MenPet-rwre,MenPetPop,LyoPem}
for instance)
that allow to obtain useful properties for the probability measure
on the words of the language in terms of ergodic properties
of the random walk.
However, the classification of these random walks 
is far from being complete and their  complete study 
remains
an interesting open problem.
Other mathematical results concerning random
walks stemming from generative grammars using more
algebraic combinatorial tools have been developped in \cite{Leroux}.

\subsection{Asynchronous Grammar-driven Processes}
A sequence of internal clocks are attached to 
subwords in  $\textsf{Dom}(\Pi)$; when they ring, at random times
distributed exponentially, the subword is transformed by a new subword
according to the allowed productions. 
In general, productions do not preserve the length of the words. If we
denote by $N_a(\alpha), a\in\mathbb{A}$ the number of letters $a$ contained
in the word $\alpha\in\mathbb{A}^n$, then the passage from generation $n$ to $n+1$
induces a multibranching process \cite{Jagers}
whose sub-populations behave like $N_a$. 
To keep evolving  words inside
the same space we consider infinite length words from the very beginning
and apply asynchronous evolution on the infinite sequence. This procedure
is quite standard in statistical physics; although mathematically more
delicate to handle, the obtained results are sharper than
the finite case;  finite size results can  be inferred from
infinite       sequences.
\begin{example}
The Fig.~\ref{fig:async} illustrates how asynchronous grammar-driven
process evolves for an infinite initial sequence.
The evolution of a fixed window of size $2N+1$ is depicted in this figure.
Notice that since productions are not length preserving
in general, there does not exist a global coordinate system to number
the residues.
\begin{figure}
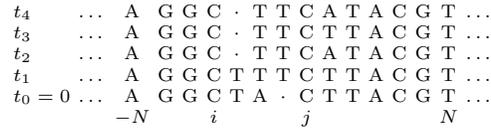
 
\begin{center}
{\scriptsize
\begin{tabular}{lcccccccccccccccc}
$t_4$ & \ldots & A&G&G&C&$\cdot$&T&T&C&A&T&A&C&G&T&\ldots\\
$t_3$ & \ldots & A&G&G&C&$\cdot$&T&T&C&T&T&A&C&G&T&\ldots\\
$t_2$ & \ldots & A&G&G&C&$\cdot$&T&T&C&A&T&A&C&G&T&\ldots\\
$t_1$ & \ldots & A&G&G&C&T&T&T&C&T&T&A&C&G&T&\ldots\\
$t_0=0$ & \ldots & A&G&G&C&T&A&$\cdot$&C&T&T&A&C&G&T&\ldots\\
    &        & $-N$&&&$i$&&&&$j$&&&&&&$N$&
\end{tabular}}
\end{center}
\caption{\label{fig:async} A random realisation of a computational path:
Productions $\alpha_1\alpha'\alpha_2\rightarrow\alpha_1\beta'\alpha_2$ 
occur
at random times.
At     $t_1$: $\alpha_1=T$, $\alpha_2=C$, $\alpha'=A$ and $\beta'=TT$.
At  $t_2$: $\alpha_1=GC$, $\alpha_2=TT$, $\alpha'=T$ and $\beta'=\kappa$.
etc. 
If we allow infinite re-numberings in order
to impose a global coordinate system, at some places the configuration
must be squeezed.
The symbol $\cdot$ reminds where squeezing takes place}
\end{figure}
\end{example}

In \cite{MalyshevRG}, stochastic evolution in the absence
of a global coordinate system has been studied and in \cite{MalyshevQG}
this method has been extended ot quantum evolution in the context of
quantum gravity. These results are  presented below adapted to a genetic
context.

\subsubsection{The Classical Stochastic Case}
 We consider infinite length configurations
in $\mathbb{X}=\mathbb{A}^\mathbb{Z}$ and the set $\Omega$ of continuous
time processes $\Omega=\{\omega:[0,\infty[\rightarrow \mathbb{X} \
\textrm{admissible}\}$. A process is termed 
\emph{admisssible} if it is right continuous and 
if $\omega(s-)\ne \omega(s)$ for some $s\in [0,\infty[$ then
there exist a left semi-infinite word $\alpha$, a right semi-infinite
word $\gamma$ and two finite words $\beta,\beta'$ with $(\beta,\beta')\in\Pi$
such that  $\omega(s-)=\alpha\beta\gamma$ and
$\omega(s)=\alpha\beta'\gamma$.
Now fix some positive integer $N$ and a finite initial
configuration $\xi$ of length
$2N+1$. In the Fig.~\ref{fig:async} above, $N=6$ and the initial configuration
$\xi=AGGCTACTTACGT$ of length 13. Notice that the symbol $\cdot$ does not
make part of the initial sequence!
Denote by $\Omega_N[\xi]=
\{\omega\in\Omega: \omega_{-N}(0)=\xi_{-N}\ldots  \omega_{N}(0)=\xi_{N}\}$
the set of process trajectories emanating from  the cylinder set
defined by the fixed configuration $\xi$.

Since a global coordinate system cannot be used, local observers
are placed, let say at the original sites $i$ and $j$ 
(see Fig.~\ref{fig:async}); denote  by $\Omega_{N; i,j;t}$
the subset of such trajectories for which the residues on sites
$i$ and $j$ have not been modified up to time $t$ while this property
is false for  sites $k$ with $i<k<j$.
The probability rates
on the productions induce the probability 
$\mathbb{P}(\Omega_{N; i,j;t})$ on the set   $\Omega_{N; i,j;t}$.

\begin{theorem}
The limit $\lim_{N\rightarrow \infty} \mathbb{P}(\Omega_{N; i,j;t})$ 
exists for all $i,j\in\mathbb{Z}, i<0\leq j$ and defines, for all
$t\in[0,\infty[$, a joint probability $\mu(i,j;t)$ such
that
\[\sum_{(i,j): i<0\leq j} \mu(i,j;t)=1.\]
\end{theorem}
\begin{remark}
Although the above theorem appears as an existence
result, as a byproduct of the proof, one gets very precise
estimates on the properties of the probability measure $\mu$.
Fixing the grammar, we get estimates
of the joint probability measure on sequences.
\end{remark}

\subsubsection{The Quantum Case}
 Productions are implemented
by operators; more precisely, suppose that 
$\alpha=\alpha_1 \alpha'\alpha_2$ with $|\alpha_1|=j-1$, for $j\geq 1$,
and $\beta=\alpha_1 \beta'\alpha_2$, while $\pi:=(\alpha',\beta')\in\Pi$.
Then we define an operator $A_\pi(j)$ by its action on basis vecors:
$A_\pi(j)e_\alpha=e_\beta$ and similarly its adjoint
$A^*_\pi(j)e_\beta=e_\alpha$.
These operators define an operator, the \emph{Hamiltonian},
\[H=\sum_{\pi\in\Pi} \sum_{j\in\mathbb{N}}(\lambda_\pi A_\pi(j)+
\overline{\lambda}_\pi A^*_\pi(j))+cH_c,\]
with $\lambda_\pi\in\mathbb{C}$, $\overline{\lambda}_\pi$
denoting the complex conjugate of $\lambda_\pi$, $c\in\mathbb{R}$, and
$H_c$ a diagonal correction term.

The above Hamiltonian is formally
 self-adjoint. Therefore the operator $U(t)=\exp(itH)$ is formally
unitary and corresponds to time evolutions.
Introducing
the family of Hilbert supspaces $\mathbb{H}_N=
\textsf{span}\{e_\alpha, |\alpha|\leq N\}$
and denoting by $P_N$ the projector
to $\mathbb{H}_N$, we get  regularised finite length
Hamiltonians $H_N=P_N H P_N$.
The set of all self-adjoint operators
on $\mathbb{H}_N$ constitutes a finite dimensional
$C^*$-algebra denoted by $\mathfrak{A}_N$. When passing to 
the inductive limit $\mathfrak{A}_\infty=
\overrightarrow{\lim}_{N\rightarrow \infty}
\mathfrak{A}_N$, and then to its norm-closure $\mathfrak{A}=\overline{\mathfrak{A}_\infty}$ 
we obtain the so-called quasi-local algebra of observables $\mathfrak{A}$.
Recall that the set of self-adjoint operators together with a state
is the quantum (non-commutative) analogue of random variables. Hence
the relevant question is whether a state can be defined on the 
 quasi-local algebra $\mathfrak{A}$. (These are 
standard constructions in the context of quantum statistical
mechanics; details can be found for instance in \cite{BratteliRobinson2}.) 

\begin{theorem}
\label{thm:qc}
Define $\rho_N=\frac{\exp(-\tau H_N)}{\mathsf{tr}\exp(-\tau H_N)}$, for
$\tau\geq 0$; then
$\rho_N$ is a state on  $\mathfrak{A}_N$. For every quantum observable
$X\in\mathfrak{A}$ we have
$\lim_{N\rightarrow\infty} \mathsf{tr}(\rho_N X)= \mathsf{tr}(\rho X)$,
where $\rho$ is a state on the infinite system.
\end{theorem}
\begin{remark} 
Giving the productions and the corresponding unitary vector of evolution
defines a state on the set of observables over infinite length sequences 
by virtue of  Theorem~\ref{thm:qc}
\end{remark}

\section{Conclusion and Open Problems}
We have presented as simply as possible a formalism based on general grammars
acting on sequence spaces. It is shown that for context-sensitive
stochastic generational grammars, the relevant object
to study is a random walk in random environment on tree. Asymptotic behaviour
for such objects is known only for very special models. It is a challenging
open problem to have a more complete classification of these random walks.

For quantum grammars, the relevant objects are quantum random walks.
Now the known models of random walks are essentially only one-dimensional.
Thus there are challenging open problems even for 
context-free and regular quantum grammars corresponding to quantum
random walks on trees.

Then we have presented the case of asynchronous (random or quantum)
grammar-driven processes on infinite sequences and stated known results
establishing the existence of a joint state on the infinite-dimensional
algebra of
observables. These results show that if we know the grammar, there
exist a global state that stems from this grammar. 

There are several
open  challenging
problems in this context. Firstly, from experimental observations
on very long sequences, we can estimate correlation properties of the
state. Is it possible to reconstruct the Hamiltionian (hence the
grammar) giving rise to this state?
A second important problem is unicity : is it true that
a given Hamiltonian gives rise to a unique state or some
phenomenon of phase transition occurs? The consequences of
a phase transition would be that the system becomes unstable;
although two different cells share the same grammar, it is enough
that some very small external perturbation acts differently on each
of them for their genetic sequence to evolve (mutate) to different states.
 
The above mentionned problems are essentially mathematical in nature.
There are however several biological and algorithmic problems associated
with them. 
The fundamental thesis defended in this work is that the genome 
statistics of a class of individuals of a given species is 
determined by the stochastic or quantum grammar inducing the 
asynchronous process. 
Assuming absence of phase transition and 
that for individuals belonging to two different classes 
(that can be distinguished for instance by an experimentally 
observed spectacular difference in the reaction to a drug) 
the genome statistics can be discriminated, 
it follows that the determining grammars must be different. 
In order to accept this thesis, the algorithmic problems 
must be solved and precisely designed biological experiments 
must confirm it. 
But if it is eventually established, it provides a 
mesoscopic explanatory scheme involving the 
fundamental mechanisms that govern the time evolution of the DNA molecule.


 \bibliographystyle{plain}
 \bibliography{petritis,rwre,bio,cours-qinfo}


\end{document}